\documentclass[oneside, 12pt]{amsart}
\usepackage{amscd, amssymb, amsmath, mathrsfs}
\usepackage[english]{babel}
\usepackage{booktabs}
\usepackage{tikz-cd}
\usepackage{url}
\usepackage[pdftex, colorlinks=true,  citecolor=blue, linkcolor=blue, linktocpage=true]{hyperref}

\newcommand\mylabel[1]{\label{#1}\marginpar{\vspace{-1ex}\medskip\medskip\footnotesize \tt #1}}
\renewcommand\mylabel[1]{\label{#1}}
\newcommand{\mydate}{
\number\day\space
\ifcase\month \or January\or February\or March\or April\or May\or June\or July\or August\or September\or October\or November\or December\fi 
\space\number\year}

\DeclareUrlCommand\arXiv{\urlstyle{same}}

\newtheorem{theorem}{Theorem}[section]

\newtheorem{proposition}[theorem]{Proposition}

\theoremstyle{definition}

\newtheorem*{acknowledgement}{Acknowledgement}

\theoremstyle{remark}

\newcommand{\PP}{\mathbb{P}}

\newcommand{\GG}{\mathbb{G}}

\newcommand{\shL}{\mathscr{L}}

\newcommand{\Aff}{\text{\rm Aff}}

\newcommand{\Aut}{\operatorname{Aut}}

\newcommand{\Bl}{\operatorname{Bl}}

\newcommand{\Gal}{\operatorname{Gal}}
\newcommand{\GL}{\operatorname{GL}}

\newcommand{\id}{{\operatorname{id}}}

\newcommand{\lra}{\longrightarrow}

\renewcommand{\O}{\mathscr{O}}

\newcommand{\op}{\text{\rm op}}

\newcommand{\Pic}{\operatorname{Pic}}

\newcommand{\quadand}{\quad\text{and}\quad}

\newcommand{\ra}{\rightarrow}

\newcommand{\sep}{{\operatorname{sep}}}
\newcommand{\Set}{{\text{\rm Set}}}

\newcommand{\Spec}{\operatorname{Spec}}

\newcommand{\uHom}{\underline{\operatorname{Hom}}}

\begin{document}

\title[Algebraic spaces that become schematic]
      {Algebraic spaces that become schematic after ground field extension}.

\author[Stefan Schr\"oer]{Stefan Schr\"oer}
\address{Mathematisches Institut, Heinrich-Heine-Universit\"at, 40204 D\"usseldorf, Germany}
\curraddr{}
\email{schroeer@math.uni-duesseldorf.de}

 
\dedicatory{Revised version, 13 September 2021}

\begin{abstract}
We construct examples of non-schematic algebraic spaces that become schemes after finite ground field extensions.
\end{abstract}

\maketitle
\tableofcontents
 
\section*{Introduction}
\mylabel{Introduction}

Roughly speaking, \emph{algebraic spaces} are generalizations of \emph{schemes} that allow more freely the  formation of contractions and   quotients,
which are  somewhat problematic in the  category of schemes.
For example, if $E$ is a connected curve on a smooth projective surface $Y$ such that the irreducible components
form a negative-definite intersection matrix $N=(E_i\cdot E_j)$, the contraction of the curve to a point
exists as an algebraic space $X$ that is usually non-schematic if $h^1(\O_E)>0$.   
Furthermore, if a finite group $G$ acts freely on a separated scheme $Y$ of finite type, the quotient
$X=G\backslash Y$ exists as an algebraic space, but it  must be  non-schematic if some orbit $G\cdot a$ does not admit an affine open neighborhood.

Despite their importance, it is by no means straightforward to provide  geometrically meaningful examples of algebraic spaces that are indeed not schemes.
The goal of this note is to highlight the ubiquity of such non-schematic algebraic spaces:
We construct algebraic spaces $X$ over non-closed  ground  fields 
$k$ that are non-schematic, but where the  base-change $X\otimes k'$ to some  finite Galois extension
become schemes.
The construction  is   easy, and relies on schemes $X_0$ on which
some finite quotient $G$ of the absolute Galois group $\Gal(k^\sep/k)$ acts such  that some orbits  do  not admit   affine open neighborhoods.
The procedure is formalized in Theorem \ref{twisting non-schematic}, which gives a criterion for Galois twists to be non-schematic.

Our first concrete example comes from a  proper normal surface $X_0$ birational to $E\times \PP^1$, 
for some elliptic curve $E$ where  the group of rational points  $E(k)$ is not torsion.
This depends on earlier constructions of myself \cite{Schroeer 1999}.

The second example is actually a general construction 
starting  from any separated normal scheme $\tilde{Y}$ of finite type that does not admit an ample invertible sheaf.
Now  $X_0$ is obtained from the  disjoint union of two copies of $\tilde{Y}$ via   an identification of certain closed points. 
The idea is extracted from a beautiful observation of  Artin (\cite{Artin 1971}, page 286). It also   relies on a result of Benoist \cite{Benoist 2013}
that there are only finitely many maximal quasiprojective open sets
in $\tilde{Y}$, which solves a conjecture by of   Bia\l{}ynicki-Birula (\cite{Bialynicki-Birula 1998}, page 302), and
generalizes earlier results of Kleiman \cite{Kleiman 1966} and  W\l{}odarczyk \cite{Wlodarczyk 1999}.

The reviewer pointed out  that     Huruguen \cite{Huruguen 2011}    studied schemes   such that some base-change  
aqcuire the structure of a toric or spherical variety, and that in this framework one may also obtain non-schematic algebraic spaces
whose base-change become toric or spherical varieties.
Let me also mention that there are families of curves $X\ra\Spec(R)$ over  discrete valuation rings $R$
where $X$ is non-schematic but the base-change to the henselization $R'$ becomes schematic
(see \cite{Bosch; Luetkebohmert; Raynaud 1990}, Section 6.7 and  \cite{Schroeer 2001}).
 
\begin{acknowledgement}
I wish to thank the reviewer for valuable comments.
This research was conducted in the framework of the   research training group
\emph{GRK 2240: Algebro-Geometric Methods in Algebra, Arithmetic and Topology}, which is funded
by the Deutsche Forschungsgemeinschaft. 
\end{acknowledgement}

\section{Non-schematic twists}
\mylabel{Non-schematic}

Let $S$ be a base scheme, and write $(\Aff/S)$ for the category of 
affine $S$-schemes.
Recall that an \emph{algebraic space} is a contravariant functor
$$
F:(\Aff/S)\lra (\Set)
$$
satisfying the sheaf axiom for the \'etale topology, such that
the diagonal monomorphism $F\ra F\times F$ is relatively representable by schemes,
and that there is an \'etale surjection $U\ra F$ from a scheme.
 Each scheme $X$ can be viewed as an algebraic space, via the Yoneda embedding,
and an algebraic space that corresponds to a scheme is called \emph{schematic}.
For more details on algebraic spaces, we refer to the monographs of Artin \cite{Artin 1973}, Knutson \cite{Knutson 1971} and  Olsson  \cite{Olsson 2016}.

A tremendous advantage of algebraic spaces is that forming quotient is permissible, in the following
very general setting:
Let $G$ is an algebraic space endowed with a group structure acting on an algebraic space $F$.
Suppose the action is free, and that the    structure morphism $G\ra S$
is flat and locally of finite type. Then the  sheaf  quotient $ G\backslash F$ is an algebraic space,
and the quotient map $F\ra G\backslash F$ is flat and locally of finite presentation 
(for example \cite{Laurent; Schroeer 2021}, Lemma 1.1).  

For the sake of exposition, we now suppose that $S$ is the spectrum of a ground field $k$.
Let $X_0$ be an algebraic space that is separated and of finite type.
Then $\Aut_{X_0/k}$ is a group scheme that is  locally of finite type.
Let $G$ be any group scheme that is locally of  finite type,    $G\ra\Aut_{X_0/k}$ be a homomorphism, and consider the resulting $G$-action on $X_0$.
For each principal homogeneous $G$-space $P$, the diagonal action on the  product $P\times X_0$ is free, hence the quotient 
$$
X={}^PX_0 = P\wedge^G X_0=G\backslash (P\times X_0)
$$
is an algebraic space. Composing the quotient map $q$ for the diagonal action with the   inclusion $\Delta_P$ of the diagonal for the principal homogeneous space, we obtain
$$
P\times X_0 \stackrel{\Delta_P\times\id_{X_0}}{\lra} (P\times P)\times X_0= P\times (P\times X_0)\stackrel{\id_P\times q}{\lra} P\times X,
$$
which is an isomorphism compatible with the projection to $P$. So for each morphism $\Spec(L)\ra P$ with some field $L$,
we get a canonical identification $X_0\otimes L= X\otimes L$. In particular, $X$ is a \emph{twisted form} of $X_0$.

We are primarily interested in the case that $G\subset\Aut(X_0)$ is a finite subgroup 
and   $P=\Spec(k')$, where $k'$ is a Galois extension with $G=\Gal(k'/k)^\op$.
Then $X={}^PX_0$ is called the \emph{Galois twist} of $X_0$ with respect to $L$. By the above, 
we have a canonical identification $X\otimes k'=X_0\otimes k'$. This leads to the following fact:
 
\begin{theorem}
\mylabel{twisting non-schematic}
In the above setting, suppose $X_0$ is schematic, and that there is a   rational point $a \in X_0$ whose orbit
$G\cdot a \subset X_0$ does not admit an affine open neighborhood.
Then the Galois twist $X$ is non-schematic,
whereas its base-change $X\otimes k'$ becomes schematic.
\end{theorem}
 
\proof
The orbit $Z_0=G\cdot a $ is finite, hence schematic.
Write  $A_0=H^0(Z_0,\O_{Z_0})$ for the coordinate ring. Then the quotient $Z=G\backslash(Z_0\otimes k')$ is the spectrum for 
the ring of invariants $L=(A_0\otimes k')^G$,
which coincides with the fixed field $L\subset k'$ for the stabilizer subgroup $H=G_a$.
This defines a morphism $\Spec(L)\ra X$. One may regard this as a point on the algebraic space, and 
we now check that $X$ is not schematic near this point.

Seeking a contradiction, we assume that there is an affine open subspace $U\subset X$ over which $\Spec(L)\ra X$ factors.
Then $U\otimes k'\subset X\otimes k'=X_0\otimes k'$ is an affine open neighborhood of $Z_0\otimes k'$.
We now use that the projection $X_0\otimes k'\ra X_0$ is the quotient with respect to the Galois action.
Clearly
$$
V'=\bigcap_{\sigma\in G} (\id_{X_0}\otimes \sigma)(U\otimes k')\subset X_0\otimes k'
$$
is a Galois-invariant affine open neighborhood of $Z_0\otimes k'$, and its quotient by the Galois action
defines an affine open neighborhood for $Z_0$ inside $X_0$, contradiction.
Thus the twisted form $X$ is non-schematic. On the other hand, the identification $X\otimes k'=X_0\otimes k'$ shows
that its base-change to $k'$ becomes schematic.
\qed

\section{Examples}
\mylabel{Examples}

Fix a ground field $k$. We now give  some geometrically meaningful examples    for which Theorem \ref{twisting non-schematic} applies.

First consider the smooth projective surface $E\times\PP^1$, where $E$ is an elliptic curve.
Fix rational points $u,v\in E$ and $\lambda,\mu\in\PP^1$.
We   recall from \cite{Schroeer 1999} how this   leads to a proper normal surface $X_0$
with   two elliptic singularities:
Write  $\Bl_Z(E\times\PP^1)$ for the blowing-up 
with respect to the center
$Z=\{( u,\lambda),(v,\mu)\}$.
Let $D_\lambda,D_\mu\subset\Bl_Z(E\times\PP^1)$ be the strict transform of the fibers $E\times\{\lambda\}$ and $E\times\{\mu\} $, respectively.
These    are copies of the elliptic curve $E$, each with self-intersection $s=-1$.
In turn,   $D=D_\lambda\cup D_\mu$ is a negative-definite curve, so 
the  contraction $f:\Bl_Z(E\times\PP^1)\ra X_0$  of $D$ exists as an algebraic space (\cite{Artin 1970}, Corollary 6.12),
which here is actually a scheme (\cite{Schroeer 1999}, Section 2.3). 
This $X_0$ is  a proper normal surface with $h^0(\O_{X_0})=1$,
and its singular locus   comprises the two rational points $a=f(D_\lambda)$ and $b=f(D_\mu)$.

\begin{proposition}
\mylabel{involution exists}
There is an involution $\sigma\in\Aut(X_0)$ with $\sigma(a)=b$.
\end{proposition}

\proof
Without loss of generality we may assume that $\lambda=(1:0)$ and $\mu=(0:1)$.
Then the matrix
$(\begin{smallmatrix}0&1\\1&0\end{smallmatrix})\in\GL_2(k)$
defines an involution $\tau_2\in\Aut(\PP^1)$ with $\tau_2(\lambda)=\mu$.
Furthermore, $x\mapsto -x+(u+v)$ is an involution $\tau_1\in\Aut(E)$ of the underlying genus-one curve with $\tau_1(u)=v$.

The diagonal involution $\tau=\tau_1\times \tau_2$ on $E\times\PP^1$ leaves the center invariant, 
thus induces an involution on $\Bl_Z(E\times\PP^1)$. This induced involution clearly permutes the fibers of the projection to 
$\PP^1$. It follows that the strict transform $D$ of   $E\times\{\lambda,\mu\}$ is 
an invariant subset, and its two connected components are permuted by the   induced involution. By the universal property of contractions,
we obtain the desired    involution $\sigma:X_0\ra X_0$ that permutes  the two singularities $a,b\in X_0$.
\qed

\medskip
Benoist (\cite{Benoist 2013}, Theorem 9) showed that  each  separated  normal  scheme  $Y$ of finite type
contains only finitely many maximal quasiprojective open sets $U_1,\ldots,U_r\subset Y$. 
This solved a conjecture  of   Bia\l{}ynicki-Birula (\cite{Bialynicki-Birula 1998}, page 302),
 and generalizes earlier results of Kleiman \cite{Kleiman 1966} and  W\l{}odarczyk \cite{Wlodarczyk 1999}.
For information  on the number $r\geq 0$, see \cite{{Farnik; Jelonek 2010}},  \cite{Farnik 2013} and \cite{Farnik 2016}.

\begin{proposition}
\mylabel{no ample sheaf}
If the difference $u-v\in E(k)$ has infinite order, then the proper normal scheme $X_0$ admits no ample invertible sheaf, and 
there are exactly $r=2$ maximal quasiprojective open sets, namely 
$$
U_1=X_0\smallsetminus\{a\}\quadand U_2=X_0\smallsetminus\{b\}.
$$
In particular, the points $a,b\in X_0$ do not admit a common affine open neighborhood.
\end{proposition}

\proof
According to \cite{Schroeer 1999}, Section 2.3 the individual contraction of $D_\lambda\subset \Bl_Z(E\times\PP^1)$ yields
a projective scheme, which contains a copy of $U_2$ as an open subscheme. Hence $U_2$ is quasiprojective.
Suppose it would not be maximal. Then $X_0$ is quasiprojective, so there is a curve $D_0\subset X_0$
disjoint from $\{a,b\}$. Its strict transform $D\subset E\times\PP^1$ is a Cartier divisor with
$$
D\cap (E\times\{\lambda\})=m_1\cdot(u,\lambda)\quadand D\cap (E\times\{\mu\})=m_2\cdot(v,\mu) 
$$
for some   $m_1,m_2\geq 1$, as explained in \cite{Schroeer 1999}, Section 2.5. Interpreting these multiplicities as intersection numbers
of the invertible sheaf $\shL=\O_{E\times\PP^1}(D)$ with   linearly equivalent curves, we see that the numbers must coincide.
Write $m\geq 1$ for the common value.
Since the canonical map  $\Pic(E)\oplus\Pic(\PP^1)\ra\Pic(E\times\PP^1)$ is bijective, one gets   $\O_E(u)^{\otimes m}\simeq \O_E(v)^{\otimes m}$.
Let $\infty\in E$ be the point at infinity, which is the zero element in the group $E(k)$.
Using the identification $E(k)=\Pic^0(E)$ given by $w\mapsto \uHom(\O_E(\infty),\O_E(w))$,
we infer that $u-v\in E(k)$ is annihilated by the integer $m\geq 1$, contradiction.
The same argument applies to $U_1$.
Summing up, the $U_1,U_2\subset X_0$ are    maximal quasiprojective open sets. 

Suppose  there is another maximal quasiprojective open set $U\neq U_i$.
It must contain both $a,b\in X_0$, hence the two points admit a common affine open neighborhood. With
\cite{Hartshorne 1970}, Chapter II, Proposition 3.1 we get a curve $D_0\subset X_0$ as above and again reach a contradiction.
\qed

\medskip
We now sketch a   completely different construction that gives   schemes $X_0$ for which Theorem \ref{twisting non-schematic} applies.
This relies on the following general procedure:  Consider any separated scheme  $\tilde{Y}$  of finite type.
Let  $Z\subset \tilde{Y}$ be a finite closed subscheme that admits an affine open neighborhood, and write    $A=H^0(Z,\O_Z)$ for the coordinate ring.
We  can form the cocartesian square
$$
\begin{CD}
\Spec(A)	@>>>	\tilde{Y}\\
@VVV		@VVqV\\
\Spec(k)	@>>a>	Y
\end{CD}
$$
in the category of sheaves on $(\Aff/k)$. Note that the upper horizontal  arrow is the closed embedding of the finite scheme  $Z= \Spec(A)=\{z_1,\ldots,z_r\}$,
and the vertical arrow to the left is  the structure morphism.
The push-out $Y$ is an algebraic space (\cite{Artin 1970}, Theorem 6.1), which  here is actually a scheme (\cite{Ferrand 2003}, Theorem 7.1). 
Such constructions are often called \emph{pinching}. By abuse of notation, we  simply say that $q:\tilde{Y}\ra Y$
is the \emph{identification that turns the  closed  points   $z_1,\ldots,z_r $  into a rational point 
$a\in Y$}. Note that $Y$ remains separated and of finite type.

Now suppose that the scheme $\tilde{Y}$ \emph{is normal but does not admit an ample invertible sheaf}.
For example, this could be the  surface in Proposition \ref{no ample sheaf},
or Hironaka's   proper smooth threefold (see \cite{Hartshorne 1977}, Appendix C, Example 3.4.1).
 Let  $U_1,\ldots,U_s\subset \tilde{Y}$, $s\geq 2$
be the maximal quasiprojective open  sets, and choose closed points $u_i\in \tilde{Y}\smallsetminus U_i$.
Clearly, the closed set $\{u_1,\ldots,u_s\}$ does not admit an affine open neighborhood.
Consequently, there must be certain closed points  $z_1,\ldots,z_r\in \tilde{Y}$ that do  not admit a common affine open neighborhood,
 now with $r\geq 2$ minimal.  The latter ensures that the $r-1$ points  $z_2,\ldots, z_r\in \tilde{Y}$ do admit such   neighborhoods.
Let $\tilde{Y}\ra Y$ be the identification that turns $z_1$ into a rational point $u\in Y$,
and also $z_2,\ldots,z_r$ into a rational point $v\in Y$. Note that this is obtained by applying the procedure in the previous paragraph twice.
Since $q:\tilde{Y}\ra Y$ is finite  and hence affine,  and $z_1,\ldots,z_r\in \tilde{Y}$ do not
admit a common affine open neighborhood, the same holds for $u,v\in Y$. In other words,
passing from $\tilde{Y}$ to $Y$ reduces from $r\geq 2$ to the case $r=2$.

Let $G= \{e,\sigma\}$ by the cyclic group of order two. Regard it as a constant group scheme,
and form $G\times Y$. This  is the disjoint union of two copies of $Y$, endowed with a permutation action of $G$.
Consider the identification $q:G\times Y\ra X_0$ 
that turns $(e,u)$ and $(\sigma,v)$ into a rational point $a\in X_0$, 
and also turns $(e,v)$ and $(\sigma,u)$ into a rational point $b\in X_0$.
For a useful illustration of the geometry of $X_0$, see \cite{Artin 1971}, Picture on page 286.

\begin{proposition}
\mylabel{action descends}
In the above situation, the rational points $a,b\in X_0$ do not admit a common affine open neighborhood,
the $G$-action on $G\times Y$ induces   a free action on $X_0$, and we have  $G\cdot a=\{a,b\}$.
\end{proposition}

\proof
Suppose there is an affine open neighborhood $U$ of $a,b\in X_0$.
Then $q^{-1}(U)$ is an affine open neighborhood of 
$$
q^{-1}(\{a,b\}) =\{(e,v),(\sigma,v), (e,u), (\sigma,u) \} = G\times\{u,v\}
$$
in $G\times Y$.
The  intersection with $\{e\}\times Y$ gives an affine open neighborhood of the points $u,v\in Y$, contradiction.
The $G$-action leaves  $q^{-1}(\{a,b\})$ invariant, while it interchanges the subsets $q^{-1}(a)$ and $q^{-1}(b)$.
With the universal property of push-outs, we infer that  the permutation action on $G\times Y$ induces a $G$-action on $X_0$,
such that  $\sigma(a)=b$.
 
Clearly, the action is free on $X_0\smallsetminus\{a,b\}$, which can be regarded as an invariant open set
of $G\times Y$. Moreover, the stabilizer groups for the rational points $a,b\in X_0$ are trivial. It follows that the $G$-action
on $X_0$ is free.
\qed

\medskip
Note that the composition   $Y=\{e\}\times Y \ra X_0\ra G\backslash X_0$ is the  identification that turns
$u,v\in Y$ into a single rational point. In particular, if $\tilde{Y}$ is   integral, the same holds for  $G\backslash X_0$.

Also note that we may have started the construction with any separated scheme $Y$ of finite type 
containing two closed points $u,v\in Y$ that do not admit a common affine open neighborhood.
For $Y$   normal and $k$   algebraically closed, this are precisely the schemes that do not admit an embedding into any toric variety,
according to   W\l{}odarczyk's  result (\cite{Wlodarczyk 1993}, Theorem A). 

On the other hand, Horrocks \cite{Horrocks 1971} showed
that the the  associated fiber bundle   $Y=E\wedge^{\GG_m} C$ stemming from a non-trivial principal $\GG_m$-bundle
$E\ra\PP^1$ and the rational nodal curve $C$ does not embed into any regular scheme, yet every finite subset admits a common affine open neighborhood.
This Cohen-Macaulay proper surface $Y$ can be seen as a Hirzebruch surface with invariant $e>0$ in which the
negative and positive section are identified.

Summing up, we have described above two constructions of separated schemes  $X_0$ of finite type,
together with an action of the cyclic group $G$ of order two,
such that some orbit $G\cdot a=\{a,b\}$ does not admit a common affine open neighborhood.
Theorem \ref{twisting non-schematic} ensures that  in both cases the twisted form  $X=(X_0\otimes k')/G$   is a  non-schematic algebraic space 
whose base-change $X\otimes k'=X_0\otimes k'$ becomes schematic, provided there is a   separable extension $k\subset k'$ of degree two.
Note that the latter indeed exists if $k$ is a finite field or a number field, or if we replace the ground field $k$ by the transcendental extension $k(T)$.


\end{document}